\documentclass{article}
\usepackage{amssymb,latexsym,amscd}

\newcommand{\se}[1]{{\section{#1}} {\setcounter{equation}{0}}}
\newtheorem{theorem}{Theorem}[section]
\newtheorem{lm}{Lemma}[section]
\newtheorem{prop}{Proposition}[section]
\newtheorem{de}{Definition}[section]

\def\k{{K\"{a}hler }}
\def\ke{{K\"{a}hler-Einstein }}

\begin{document}
\hbadness=10000
\title{{\bf Degeneration of \ke hypersurfaces in complex torus to generalized pair of pants decomposition}}
\author{Wei-Dong Ruan\\
Department of Mathematics\\
University of Illinois at Chicago\\
Chicago, IL 60607\\}
\footnotetext{Partially supported by DMS-0104150.}
\maketitle
\se{Introduction}
According to the conjecture of Calabi, on a complex manifold $X$ with ample canonical bundle $K_X$, there should exist a \ke metric $g$. Namely, a metric satisfying ${\rm Ric}_g = -\omega_g$, where $\omega_g$ is the \k form of the \k metric $g$. The existence of such metric when $X$ is compact was proved by Aubin and Yau (\cite{Yau2}) using complex Monge-Amp\`{e}re equation. This important result has many applications in \k geometry. Starting with this important result, Yau initiated the program of applying \ke metrics to algebraic geometry (\cite{yau}). It was realized by him the need to study such metrics for quasi-projective manifolds (\cite{yau1}) and their degenerations. The original proof of \cite{Yau2} was a purely existence result. Later the existence of $C_1(X)<0$ \ke metrics was generalized to complete complex manifolds by Cheng and Yau (\cite{CY}), where other than existence the proof also exhibits the asymptotic behavior of the \ke metric near the infinity boundary.\\

Since \ke metric is canonical for a complex manifold, one would expect its structure to be closely related to the topology and complex geometry of the manifold. In this work we will explore one such relation between the convergence of \ke manifolds in the sense of Cheeger-Gromov and the algebraic degeneration of the underlying algebraic manifolds.\\

From algebraic geometry point of view, when discussing the compactification of the moduli space of complex manifold $X$ with ample canonical bundle $K_X$, it is necessary to consider a holomorphic degeneration family $\pi: {\cal X} \rightarrow B$ such that $X_t = \pi^{-1}(t)$ are smooth except for $t=0$, and such that the canonical bundle of $X_t$ for $t\not=0$ as well as the dualizing sheaf of $X_0$ are ample line bundles. Such degeneration will be called canonical degeneration.\\

For $t\not=0$, the works of Aubin and Yau imply the existence of a unique \ke metric $g_t$ on $X_t$. For $t=0$, assume the central fibre $\displaystyle X_0 = \bigcup_{i=1}^l D_i$. Under certain conditions, the work of Cheng and Yau implies the existence of a complete \ke metric $g_{0,i}$ on each $\displaystyle D_i\setminus \bigcup_{j\not=i} D_j$. It is interesting to understand the relation between $g_t$ for $t\not=0$ and the collection of complete $g_{0,i}$'s. Such understanding will provide structure results on the \ke metric $g_t$ (which was only known to exist previously) based on the structure of $g_{0,i}$'s. G. Tian made the first important contribution (in \cite{Tian1}) in this direction. He proved (in \cite{Tian1}) that the family of \ke metrics $g_t$ on $X_t$ is convergent in the sense of Cheeger-Gromov to the complete \ke metric $g_0 = \{g_{0,i}\}_{i=1}^l$ on the smooth part of $X_0$ under the following three assumptions: \\

(1). the total space ${\cal X}$ is smooth;\\
(2). the central fibre $\displaystyle X_0 = \bigcup_{i=1}^l D_i$ has only normal crossing singularities; \\
(3). any three of the $D_i$'s have empty intersection.\\

Following Tian's framework authors of \cite{Lu} and later \cite{ke1} improved Tian's result by removing assumption (3). Notice that for most natural examples of canonical degeneration, the total space ${\cal X}$ is not smooth, and there are a lot of canonical degeneration singularities that cannot be reduced to normal crossing case. Therefore, there are very few natural examples beyond the curve case where the normal crossing results (\cite{Tian1,Lu,ke1}) apply. \\

One of the main purposes of this work is to give a class of interesting concrete global examples, where the theory of the degeneration of \ke manifolds applies. More specifically, we will discuss the convergence of complete \ke hypersurfaces in complex torus in the sense of Cheeger-Gromov. In the situations we are interested, the central fiber $X_0$ of the underlying algebraic degeneration is almost never normal crossing globally. In fact, the underlying algebraic degenerations are naturally toroidal canonical degenerations. In \cite{ke2}, we generalized result and method in \cite{ke1} to the case of toroidal degeneration. Methods developed in \cite{ke2} will enable us to deal with degeneration discussed in this paper. Notice that the large class of natural examples discussed in this paper concern the convergence of complete \ke manifolds with finite volume. In \cite{kesl}, we will discuss another large class of natural examples concerning the convergence of complete \ke manifolds with infinite volume.\\

The examples discussed in this paper are inspired by our previous work \cite{N} on the ``large complex structure limit" degeneration phenomenon that was first discovered by Viro \cite{V} in his work of curve patching. This original idea of Viro has many interesting applications in very diverse problems. For example, it was applied by Mikhalkin \cite{M} to real algebraic curves, by Sturmfels \cite{S} to solving algebraic equations and by our previous work \cite{N} to Lagrangian torus fibrations and mirror symmetry. From (symplectic) topology point of view, \cite{N} and \cite{M1} indicate that complex hypersurfaces in complex torus admit ``generalized pair of pants" decomposition induced by such degeneration in terms of singular (Lagrangian) fibration. For example, $C_t = \{(z_1,z_2)\in (\mathbb{C}^*)^2: p_t(z_1,z_2)=0\}$ with
\[
p_t(z_1,z_2) = t^3(1+z_1^3+z_2^3) + t(z_1 + z_1^2 + z_2 + z_2^2 + z_1z_2^2 + z_1^2z_2) + z_1z_2
\]
defines a family of cubic curves in complex 2-torus, which topologically are elliptic curves with 9 points removed. When $t$ goes to zero, $C_t$ degenerate into 9 ``pair of pants" ($\mathbb{CP}^1$ with 3 points removed that can be identified with $\{z_1+z_2+1\} \subset (\mathbb{C}^*)^2$). This example can also be viewed as an explicit realization of Deligne-Mumford stable degeneration of elliptic curves with 9 marked points into a stable curve as the union of 9 pair of pants. (For more detail of this degeneration construction and its higher dimensional generalization, please refer to section 3.) Our current work (theorem \ref{dc}) will show that the convergence of the corresponding canonical \ke hypersurfaces in the sense of Cheeger-Gromov will canonically degenerate the underlying manifolds into ``generalized pair of pants" decomposition. This result brings yet another interesting relation between the \ke metric and the complex, symplectic geometry and topology of the underlying manifold.\\

Another main objective of this paper is to find interesting examples of minimal Lagrangian submanifolds in \ke manifolds with negative first Chern class. Such minimal Lagrangian submanifolds were first discussed in \cite{B}. Up to now, most known examples of such minimal Lagrangian submanifolds are either the fixed point set of an anti-holomorphic automorphism of the \ke manifold or explicit local examples. In our opinion, such minimal Lagrangian submanifolds arise most naturally as the vanishing cycles of degeneration of \ke manifolds with negative first Chern class. In \cite{hl}, we developed deformation techniques to construct such minimal Lagrangian vanishing cycles in general. Applying techniques from \cite{hl}, we are able to represent the vanishing cycles of the degeneration discussed in this paper by minimal Lagrangian tori in the corresponding \ke hypersurfaces (theorem \ref{fa}). In the curve case, the minimal Lagrangian tori are just minimal geodesic circles that will canonically divide smooth $C_t$ near the degeneration into a union of hyperbolic pairs of pants.\\

The paper is organized as follows. Section 2 summarizes the basic facts from Riemannian geometry that are needed for our work. Section 3 introduces the algebraic framework and structure of the degenerating family $\{X_t\}$ of hypersurfaces and discusses its canonical limit $\tilde{X}_0$. In sections 4,5,6, we carry out the proof of the convergence of \ke hypersurfaces in complex torus in the sense of Cheeger-Gromov to generalized pair of pants decomposition. Since the degenerating hypersurfaces in this paper are not compact, results in \cite{ke2} do not apply directly. In section 4, proposition \ref{bb} is proved based on lemma \ref{bd}, which is the main technical result that enables us to generalize estimates in \cite{ke2} to the complete manifold case here. In section 5, the approximate metrics are constructed rather explicitly, which is an improvement over the formula in \cite{ke2} that would not apply to the complete manifold case. The construction and estimates of the approximate metrics here are global in nature and are simpler than the counterparts in \cite{ke2}. Section 6 is essentially the same as its counterpart in \cite{ke2}. In section 7, we clarify the relation between the ``large complex structure limit" degeneration in this paper and general toroidal degeneration discussed in \cite{ke2}. In section 8, we apply general result from \cite{hl} to construct minimal Lagrangian vanishing torus in our \ke hypersurfaces here. \\ 

{\bf Convention of notations:} (1) When we use $g_t$ to denote a \k metric, we will automatically use $\omega_t$ to denote the corresponding \k form and vice versa. (2) When the \k potential is $\log[V_t/(\Omega_t\wedge\bar{\Omega_t})]$ for volume $V_t$ and some holomorphic section $\Omega_t$ of $K_{X_t}$, since the \k metric is independent of the choice of $\Omega_t$, by slight abuse of terminology, we will call the \k potential the logarithm of the volume form $V_t$ and denote by $\log V_t$. (3) By $A\sim B$, we mean that there exist constants $C_2>C_1>0$ such that $C_1B \leq A \leq C_2B$. (4) Following convention in analysis, $C$ is used to denote a constant that may differ in different formulas. (5) A smooth function $f(x)$ is called a bounded smooth function of $x$ if $f(x)$ and all its multi-derivatives are bounded when $x$ is bounded. We have the following basic properties of bounded smooth functions:\\
\begin{prop}
\label{aa}
(i) If both $f$ and $g$ are bounded smooth functions, then $f\circ g$ is also a bounded smooth function. (ii) If $f(x)$ is a bounded smooth function of $x$ and $|f(x)|>C>0$, then $\log f(x)$ is a bounded smooth function of $x$. (Consequently, $1/f(x) = \exp(-\log f(x))$ is a bounded smooth function of $x$.)
\hfill\rule{2.1mm}{2.1mm}
\end{prop}

\se{Background from Riemannian geometry}
Notations in this section will not extend to other sections of this paper. Results in this section are all wellknown basic facts from Riemannian geometry. We present them here due to our failure to find a source that is written in a convenient enough form for us to quote from.\\

Let $(X,g)$ be a Riemannian $n$-manifold, and $\{v_i\}_{i=1}^n$ be vector fields that form a frame field. Let $\alpha_i$ be the dual 1-form fields, and $g = g_{ij}\alpha_i\alpha_j$. $[v_i,v_j] = a_{ij}^kv_k$. Let $\{w_i\}_{i=1}^n$ be another frame field and $w_i = b_i^jv_j$.\\
\begin{de}
The frame field $\{v_i\}_{i=1}^n$ is called proper (with respect to $g$) if $a_{ij}^k$, $g_{ij}$ and all their multi-derivatives with respect to $\{v_i\}_{i=1}^n$ are bounded and $\det(g_{ij})\geq C>0$. We say that $\{v_i\}_{i=1}^n$ is equivalent to $\{w_i\}_{i=1}^n$ if $b_i^j$ and all their multi-derivatives with respect to $\{v_i\}_{i=1}^n$ are bounded and $\det(b_i^j)\geq C>0$.\\
\end{de}
It is easy to verify the following

\begin{prop}
\label{gc}
Two proper frame fields $\{v_i\}_{i=1}^n$ and $\{w_i\}_{i=1}^n$ are equivalent to each other. A frame field $\{w_i\}_{i=1}^n$ that is equivalent to a proper frame field $\{v_i\}_{i=1}^n$ is also proper.
\hfill\rule{2.1mm}{2.1mm}\\
\end{prop}
A Riemannian manifold $(X,g)$ is said to have $C^\infty$-bounded curvature if the curvature of $(X,g)$ and all its covariant multi-derivatives are bounded. A smooth tensor $T$ on a Riemannian manifold $(X,g)$ with $C^\infty$-bounded curvature is called $C^\infty$-bounded with respect to $(X,g)$ (or $C^\infty_g$-bounded) if $T$ and all its covariant multi-derivatives are bounded. We have\\
\begin{prop}
\label{ga}
$(X,g)$ has $C^\infty$-bounded curvature if for any point $x\in X$, there exists a proper frame field $\{v_i\}_{i=1}^n$ on a neighborhood $U_x$ of $x$. A smooth tensor $T$ on such $(X,g)$ is $C^\infty$-bounded if and only if $T$ and all its multi-derivatives with respect to the proper frame $\{v_i\}_{i=1}^n$ are bounded in each $U_x$.
\end{prop}
{\bf Proof:} The connection coefficients $\Gamma_{ij}^k$ can be expressed as rational functions of $a_{ij}^k$, $g_{ij}$ and their multi-derivatives with respect to $\{v_i\}_{i=1}^n$, where denominators can only be powers of $\det(g_{ij})$. Consequently, the curvature of $(X,g)$ and all its covariant multi-derivatives can be expressed as rational functions of $a_{ij}^k$, $g_{ij}$ and their multi-derivatives with respect to $\{v_i\}_{i=1}^n$, where denominators can only be powers of $\det(g_{ij})$. Therefore, when $\{v_i\}_{i=1}^n$ is a proper frame field on $U_x$, $(X,g)$ has $C^\infty$-bounded curvature on $U_x$.\\

Notice that the difference of a covariant multi-derivative of $T$ and the corresponding multi-derivative of $T$ with respect to $\{v_i\}_{i=1}^n$ will be a bilinear combination of (connection coefficients $\Gamma_{ij}^k$ and their multi-derivative with respect to $\{v_i\}_{i=1}^n$) and (lower order multi-derivative of $T$ with respect to $\{v_i\}_{i=1}^n$). By induction, we get the second part of the proposition.
\hfill\rule{2.1mm}{2.1mm}\\

$(X,g)$ is said to be $C^\infty$-quasi-isometric to $(X,g')$ if $(X,g)$ is quasi-isometric to $(X,g')$, $g$ has $C^\infty$-bounded curvature and $g'$ is $C^\infty$-bounded with respect to $g$.\\
\begin{prop}
\label{gb}
$(X,g)$ is $C^\infty$-quasi-isometric to $(X,g')$ if and only if for any point $x\in X$, there exists a frame field $\{v_i\}_{i=1}^n$ on a neighborhood $U_x$ of $x$ that is proper with respect to both $g$ and $g'$. A smooth tensor $T$ on $X$ is $C^\infty$-bounded with respect to $g$ if and only if $T$ is $C^\infty$-bounded with respect to $g'$.
\end{prop}
{\bf Proof:} Assume that $(X,g)$ is $C^\infty$-quasi-isometric to $(X,g')$ and $\{v_i\}_{i=1}^n$ is a proper frame field with respect to $g$. Then $\det(g'_{ij}) \geq C_1\det(g_{ij}) \geq C_2 >0$ and $g'$ is $C^\infty$-bounded with respect to $g$. By proposition \ref{ga}, $g'_{ij}$ and all their multi-derivatives with respect to $\{v_i\}_{i=1}^n$ are bounded. Hence $\{v_i\}_{i=1}^n$ is a proper frame field with respect to $g'$.\\

Conversely, if $\{v_i\}_{i=1}^n$ is proper with respect to both $g$ and $g'$, then $g$ is quasi-isometric to $g'$ and both $g$ and $g'$ have $C^\infty$-bounded curvature (proposition \ref{ga}). Further more, by proposition \ref{ga}, $g'$ being $C^\infty_{g'}$-bounded implies that $g'$ and all its multi-derivatives with respect to the proper frame $\{v_i\}_{i=1}^n$ are bounded with respect to $\{v_i\}_{i=1}^n$. Consequently, $g'$ is also $C^\infty_g$-bounded. Therefore, $(X,g)$ is $C^\infty$-quasi-isometric to $(X,g')$.\\  

By proposition \ref{ga}, $T$ on $X$ is $C^\infty$-bounded with respect to $g$ if and only if $T$ and all its multi-derivatives with respect to $\{v_i\}_{i=1}^n$ are bounded if and only if $T$ is $C^\infty$-bounded with respect to $g'$.
\hfill\rule{2.1mm}{2.1mm}\\
\begin{prop}
\label{ge}
$C^\infty$-quasi-isometry is an equivalence relation.
\end{prop}
{\bf Proof:} The first part of proposition \ref{gb} implies that $C^\infty$-quasi-isometry is a symmetric relation. Assume that $(X,g)$ is $C^\infty$-quasi-isometric to $(X,g')$ and $(X,g')$ is $C^\infty$-quasi-isometric to $(X,g'')$, then $g''$ is $C^\infty$-bounded with respect to $g'$. By the second part of proposition \ref{gb}, we have that $g''$ is also $C^\infty$-bounded with respect to $g$. Consequently, $(X,g)$ is $C^\infty$-quasi-isometric to $(X,g'')$.
\hfill\rule{2.1mm}{2.1mm}\\
\begin{prop}
\label{gd}
Assume that $(X,g_t)$ are $C^\infty$-quasi-isometric to each other uniformly for different $t$ near $0$, and $\displaystyle \lim_{t\rightarrow 0} g_t(x) = g_0(x)$ for any $x\in X$. Then $g_t$ converge to $g_0$ on $X$ in $C^k$-topology on the space of Riemannian metrics as $t$ goes to $0$ for any $k$.
\end{prop}
{\bf Proof:} Since $g_t$ are $C^\infty$-quasi-isometric to each other uniformly for different $t$. For any sequence $g_{t_k}$ with $\displaystyle\lim_{k\rightarrow \infty}t_k =0$, there exists a subsequence $g_{t_{k_i}}$ that converges (necessarily to $g_0$ by the assumption of the proposition) in $C^k$-topology on the space of Riemannian metrics as $i$ goes to $\infty$ for any $k$. Consequently, $g_t$ converge to $g_0$ on $F$ in $C^k$-topology on the space of Riemannian metrics as $t$ goes to $0$ for any $k$.
\hfill\rule{2.1mm}{2.1mm}\\

\se{Basic setting}
Consider an integral convex polyhedron $\Delta$ (in a rank $l$ lattice $M$) with a real valued convex function $w = \{w_m\}_{m\in \Delta}$ that determines a simplicial decomposition of $\Delta$. Let $Z$ ($Z_{\rm top}$) denote the set of (top dimensional) simplices. Such $w$ is clearly generic. We have the family of complex hypersurface $X_t = \{s_t^{-1}(0)\} \subset N_{\mathbb{C}^*} \cong (\mathbb{C}^*)^l$, where $N = M^\vee$, $\displaystyle s_t = \sum_{m\in \Delta} t^{w_m}s_m$, $s_m$ denotes the monomial on $N_{\mathbb{C}^*} \cong (\mathbb{C}^*)^l$ corresponding to $m \in M$.\\

Let $\Omega$ be the canonical holomorphic volume form on $N_{\mathbb{C}^*} \cong (\mathbb{C}^*)^l$. For $t\not=0$ and $m\in \Delta$, we have $s_t/s_m: N_{\mathbb{C}^*} \rightarrow \mathbb{C}$. $\Omega_{t,m} = (t^{w_m}\Omega \otimes (d(s_t/s_m))^{-1})|_{X_t}$ defines a section of $K_{X_t} \cong (K_{N_{\mathbb{C}^*}}\otimes K_{\mathbb{C}}^{-1})|_{X_t}$. $\{\Omega_{t,m}\}_{m\in \Delta}$ are holomorphic sections of $K_{X_t}$ with at most logarithmic singularities at the infinity.

\[
e_t = \{\Omega_{t,m}\}_{m\in \Delta}: X_t \rightarrow \mathbb{CP}^{|\Delta|-1}
\]

is an embedding that equals to the restriction to $X_t$ of the natural embedding 

\[
i_t = \{t^{w_m}s_m\}_{m\in \Delta}: (\mathbb{C}^*)^l \rightarrow \mathbb{CP}^{|\Delta|-1}.
\]

Let $Y_t = {\rm Image}(i_t)$ and $H = \{\tilde{z}\in \mathbb{CP}^{|\Delta|-1}:\sum_{m\in \Delta} \tilde{z}_m = 0\}$. Then $X_t \cong Y_t \cap H$. $Y_t$ has a natural set-theoretical limit $Y_0$ in $\mathbb{CP}^{|\Delta|-1}$. 

\begin{equation}
Y_0 = \bigcup_{S\in Z} T_S,\ T_S = \{\tilde{z}\in \mathbb{CP}^{|\Delta|-1}:\tilde{z}_m \not= 0\ (m\in S);\ \tilde{z}_m = 0\ (m\not\in S)\}.
\end{equation}

Let $X_0 = Y_0\cap H$. We have $\displaystyle X_0\setminus{\rm Sing}(X_0) = \bigcup_{S\in Z_{\rm top}} X_{0,S}$, where $X_{0,S} = T_S\cap H$. As set-theoretical limit of $X_t$, $X_0$ generally has multiplicities. \\

To get the canonical multiplicity 1 algebraic limit of the family $\{X_t\}$ when $t \rightarrow 0$, we will need the following basic construction of an algebraic variety $Y_\Sigma$ determined by a fan $\Sigma\in M$. For each $\sigma\in \Sigma$, there is an affine variety $A_\sigma = {\rm Spec}(\mathbb{C}[\sigma])$. For $\sigma,\sigma'\in \Sigma$ satisfying $\sigma\subset \sigma'$, there is a natural semi-group morphism $\sigma' \rightarrow \sigma$ that restricts to identity map on $\sigma \subset \sigma'$ and restricts to zero map on $\sigma'\setminus \sigma$, which induce the map $h_{\sigma\sigma'}: A_\sigma \rightarrow A_{\sigma'}$. Using $\{h_{\sigma\sigma'}\}_{\sigma,\sigma'\in \Sigma}$, we may glue the affine pieces $\{A_{\sigma}\}_{\sigma\in \Sigma}$ into the singular variety $Y_\Sigma$. We have the following natural canonical (Whitney) stratification

\begin{equation}
Y_\Sigma = \bigcup_{\sigma\in \Sigma} T_\sigma,\ {\rm where}\ T_\sigma = (N/\sigma^\perp)\otimes_\mathbb{Z}\mathbb{C}^*.
\end{equation}

For each polyhedron $S\in Z$ and a vertex $m\in S$, there is the natural tangent cone $\sigma_{m,S} \subset M$ of $S$ at $m$. For each $m\in \Delta$, $\Sigma_m = \{\sigma_{m,S}\}_{m\in S\in Z}$ forms a fan in $M$. The polyhedron $S$ determines a toric variety $P_S$. The fan $\Sigma_m$ determine a singular variety $Y_{\Sigma_m}$. For $m\in S$, the affine variety $A_{\sigma_{m,S}}$ has natural embeddings $A_{\sigma_{m,S}} \hookrightarrow Y_{\Sigma_m}$ and $A_{\sigma_{m,S}} \hookrightarrow P_S$. Using such embeddings, we may glue $\{Y_{\Sigma_m}\}_{m\in \Delta}$ and $\{P_S\}_{S\in Z}$ together to form a variety $Y_Z$ such that $Y_{\Sigma_m}\cap P_S = A_{\sigma_{m,S}}$. We have the following natural canonical (Whitney) stratification

\begin{equation}
Y_Z = \bigcup_{S\in Z} \tilde{T}_S,\ {\rm where}\ \tilde{T}_S = N_S\otimes_\mathbb{Z}\mathbb{C}^*,\ N_S = N/S^\perp.
\end{equation}

For each $S\in Z$, there is a natural finite cover $i_{0,S} = \{s_m\}_{m\in S}: \tilde{T}_S \rightarrow T_S \subset \mathbb{CP}^{|\Delta|-1}$. Together, we have the global map $i_0 = \{i_{0,S}\}_{S\in Z}: Y_\Sigma \rightarrow Y_0 \subset \mathbb{CP}^{|\Delta|-1}$. Let $\tilde{X}_0 = i_0^{-1}(X_0) = \bigcup_{S\in Z} \tilde{X}_{0,S}$, where $\tilde{X}_{0,S} = i_{0,S}^{-1}(X_{0,S})$. We have $e_{0,S} = i_{0,S}|_{\tilde{X}_{0,S}}: \tilde{X}_{0,S} \rightarrow X_{0,S} \subset \mathbb{CP}^{|\Delta|-1}$ and $e_0 = i_0|_{\tilde{X}_0}: \tilde{X}_0 \rightarrow X_0 \subset \mathbb{CP}^{|\Delta|-1}$. For each $S\in Z_{\rm top}$, $\tilde{T}_S \cong (\mathbb{C}^*)^l$, under the toric gauge $w_m=0$ for $m\in S$, we have $\displaystyle i_{0,S} = \lim_{t\rightarrow 0} i_t$ and $\displaystyle e_{0,S} = \lim_{t\rightarrow 0} e_t$. In such sense, we may think of $\tilde{X}_0$ (resp. $Y_Z$) as the canonical multiplicity 1 algebraic limit ({\bf canonical limit} for short) of the family $\{X_t\}$ (resp. $\{Y_t\}$) when $t \rightarrow 0$. Each $\tilde{X}_{0,S}$ for $S\in Z_{\rm top}$ is a so-called ``generalized pair of pants", which is a finite abelian cover of $\{z\in (\mathbb{C}^*)^l|z_1 + \cdots + z_l +1 =0\}$.\\

{\bf Remark:} It is easy to see that all constructions in this section can also be carried out (with slight modification), when $w = \{w_m\}_{m\in \Delta}$ is not generic, namely, $Z$ is a convex polyhedron decomposition instead of a simplicial decomposition for $\Delta$.\\

\se{Basic estimates}
In the estimates of the later sections, we quite often need to prove certain functions on $X_t$ are $C^\infty$-bounded (uniformly with respect to $t$). The main goal of this section is to prove proposition \ref{bb}, which provides the technique for such purpose. The proof of proposition \ref{bb} depends on the estimates (lemmas \ref{bc}, \ref{bd}) for the convex function $w$.\\

For $m\in \Delta$, define
\[
a_m = \kappa -\log \eta_m,\ {\rm where}\ \eta_m = \|t^{w_m}s_m\|_t^2,\ \|s\|_t^2 = |s|^2\left(\sum_{m\in \Delta} |t|^{2w_m} |s_m|^2\right)^{-1},
\]

and $\kappa>0$ is a constant that will be determined later to make $a_m$ suitably large.\\
\begin{lm}
\label{bc}
There exists a constant $a>0$ such that for any $x\in N_{\mathbb{C}^*}$ the set

\[
S_{x,a} = \{m\in \Delta | \eta_m(x) > t^a\}
\]

is a simplex in $Z$.
\end{lm}

{\bf Proof:} According to the definition of $\eta_m$, it is easy to see that $\tilde{w} = \{\tilde{w}_m\}_{m\in \Delta}$ equals to the piecewise linear convex function $w = \{w_m\}_{m\in \Delta}$ up to the adjustment of an affine function, where $\tilde{w}_m = (\log \eta_m) /(2\log t)\geq 0$. Assume there is a subset $\tilde{S} \subset S_{x,a}$ that forms a simplex not in $Z$. Since $0 \leq \tilde{w}_m \leq a$ for $m \in \tilde{S}$, we may adjust $\tilde{w}$ by an affine function so that $\tilde{w}_m =0$ for $m \in \tilde{S}$ and $\tilde{w}_m \geq -C_1a$ for $m \not\in \tilde{S}$.\\

Generally, there exists a constant $C_2>0$ that only depends on the equivalence class of the strictly convex $w$ modulo affine functions, such that $\inf_{m\in \Delta} w_m \leq -C_2$ for any simplex $S\not\in Z$ and adjustment of $w$ by affine function (still denote by $w$) satisfying $w_m =0$ for $m \in S$. If we take $a < C_2/C_1$ and $S = \tilde{S} \not\in Z$, we have a contradiction. Therefore $S_{x,a}$ is a simplex in $Z$.
\hfill\rule{2.1mm}{2.1mm}\\
\begin{lm}
\label{bd}
For any $x\in N_{\mathbb{C}^*}$, there exists $S_x \in Z_{\rm top}$ with the filtration $S_1(=S_{x,a})\subset \cdots \subset S_K (= S_x)$, numbers $t_1(=t^a)\geq \cdots \geq t_K >0$, and $b_k,c_k>0$ for $1\leq k \leq K$ such that $\eta_m(x) \geq t_k^{c_k}$ for $m\in S_k$ and $t_k^{b_k} \geq \eta_m(x)$ for $m\in (\Delta\setminus S_k) \cap (M_k\setminus M_{k-1})$, where $M_k = {\rm Span}(S_k)$, $M_0 = \{0\}$.
\end{lm}
{\bf Proof:} Take $S_1 = S_{x,a}$ and $t_1 = t^a$, lemma \ref{bc} implies that $\eta_m(x) \geq t_1$ for $m\in S_1$ and $t_1 \geq \eta_m(x)$ for $m\in \Delta\setminus S_1$. Assume the lemma is true up to $k$. Assume that $\eta_m(x)$ reaches maximum $t_{k+1}$ at $m=m_{k+1}$ for $\displaystyle m \in \bigcup_{S\in Z,S_k\subset S} S\setminus S_k$. Let $S_{k+1} = S_k\cup \{m_{k+1}\}$ and $M_{k+1} = {\rm Span}(S_{k+1})$. Then it is easy to see that there exist $b_{k+1},c_{k+1}>0$ such that $\eta_m(x) \geq t_{k+1}^{c_{k+1}}$ for $m\in S_{k+1}$ and $t_k^{b_{k+1}} \geq \eta_m(x)$ for $m\in \Delta \cap (M_{k+1}\setminus M_k)$. By induction, we get the desired filtration and $S_x$.
\hfill\rule{2.1mm}{2.1mm}\\

Let $S = \{m_0,m_1,\cdots,m_l\} \in Z_{\rm top}$. Without loss of generality, we may normalize $w$ so that $w_{m_i}=0$ for $0\leq i \leq l$. Then by convexity of $w$, we have $w_m>0$ for $m\not\in S$. Without loss of generality, we may assume that $a_{m_i}$ is in ascending order. Take coordinate $z = (z_1,\cdots,z_l)$, where $z_i = s_{m_i}/s_{m_0}$. Then $1 \geq|z_1|\geq\cdots \geq |z_l|$. We may identify $m_0$ to be the origin of $M$. Then

\[
s_m/s_0 = z^m = \prod_{i=1}^l z_i^{m^i},\ \ {\rm for}\ m\not\in S,
\]

where $(m^1, \cdots, m^l)$ is the coordinate of $m$ with respect to the basis $\{m_1,\cdots,m_l\}$ of $M$. Let $I_m = \{i|1\leq i \leq l,\ m^i\not=0\}$ and $W^0_j = a_{m_j}z_j\displaystyle\frac{\partial}{\partial z_j}$.\\ 
\begin{lm}
\label{be}
Each derivative of a term in the following 
\begin{equation}
\label{bf}
\begin{array}{l}
\displaystyle \frac{1}{a_{m_j}},\ z_jP(a_{m_j}),\ \bar{z}_jP(a_{m_j}),\ \ \ {\rm for}\ 1\leq j \leq l;\\t^{w_m}z^mP(\{a_{m_j}\}_{j\in I_m}),\ \bar{t}^{w_m}\bar{z}^mP(\{a_{m_j}\}_{j\in I_m}),\ \ {\rm for}\ m\in \Delta\setminus S;\\\displaystyle \frac{a_{m_j}}{a_m},\ \ \ {\rm for}\ m\in \Delta\setminus S\ {\rm and}\ j \in I_m;
\end{array}
\end{equation}
with respect to $\{W^0_j,\bar{W}^0_j\}_{j=1}^l$ is a finite sum of products of a term of the same form and a bounded smooth function of terms in (\ref{bf}). (Here $P$ denote polynomials. A term of the same form as $\frac{a_{m_j}}{a_m}$ will just be the same term, while a term of the same form as $z_jP(a_{m_j})$ may have different polynomial $P$.)
\end{lm}
{\bf Proof:} 
\[
W^0_k(\log z_j) = \delta_{jk}a_{m_j}.
\]
\begin{equation}
\label{ba}
W^0_k(\log a_{m_j}) = \delta_{jk} - (a_{m_k}|z_k|^2)\frac{1}{a_{m_j}}\frac{1}{1+|z|^2}.
\end{equation}
\[
W^0_k(\log a_m) = m^k\frac{a_{m_k}}{a_m} - (a_{m_k}|z_k|^2)\frac{1}{a_m}\frac{1}{1+|z|^2}.
\]

The conclusion of the lemma is an easy consequence of these computations.
\hfill\rule{2.1mm}{2.1mm}\\
\begin{lm}
\label{bg}
There exists a constant $C>0$ such that when restricted to $X_t$, near $x\in X_t$ satisfying $S_x=S$, we have $1\geq |z_1| \geq C$.
\end{lm}
{\bf Proof:} When restricted to $X_t$,

\[
1+ z_1 + \cdots + z_l + \sum_{m\in \Delta\setminus S} t^{w_m}z^m = s/s_0 =0.
\]

Since we assumed that $|z_1|$ is the largest among $|z_i|$ for $1\leq i \leq l$, and when $S_x=S$, $m\in \Delta\setminus S$ is not in $S_{x,a}\subset S_x$. Consequently, $|t^{w_m}z^m| = |s_m/s_0| \leq |s_{m_1}/s_0| =|z_1|$ for $m\in \Delta\setminus S$. Therefore

\[
1 = -\left(z_1 + \cdots + z_l + \sum_{m\in \Delta\setminus S} t^{w_m}z^m\right) \leq |\Delta||z_1|,\ \ 1\geq |z_1| \geq \frac{1}{|\Delta|}.
\]
\hfill\rule{2.1mm}{2.1mm}\\

When restricted to $X_t$, under the coordinate $(z_2,\cdots,z_l)$, let $W_j = a_{m_j}z_j\displaystyle\frac{\partial}{\partial z_j}$ for $2\leq j \leq l$.
\begin{lm}
\label{bh}
When restricted to $X_t$, near $x\in X_t$ satisfying $S_x=S$, each derivative of a term in (\ref{bf}) with respect to $\{W_j,\bar{W}_j\}_{j=2}^l$ is a finite sum of products of a term of the same form and a bounded smooth function of terms in (\ref{bf}).
\end{lm}
{\bf Proof:} When restricted to $X_t$, near $x\in X_t$ satisfying $S_x=S$, 

\[
z_1 = -1- z_2 - \cdots - z_l - \sum_{m\in \Delta\setminus S} t^{w_m}z^m
\]

is clearly a bounded smooth function of terms in (\ref{bf}). Since $1\geq |z_1| \geq C>0$ (lemma \ref{bg}), by proposition \ref{aa}, $z_1$, $\frac{1}{z_1}$, $\log z_1$ and their complex conjugates are all bounded smooth functions of terms in (\ref{bf}). Let $f$ be a bounded smooth function of terms in (\ref{bf}), then

\[
W_k(f) = W^0_k(f) + W_k(z_1)\frac{\partial f}{\partial z_1}.
\]

The only term in $W_k(f)$ that needs comment is

\[
W_k(z_1)\frac{\partial f}{\partial z_1} = a_{m_k}z_k\frac{\partial z_1}{\partial z_k}\frac{\partial f}{\partial z_1} = -\frac{1}{a_{m_1}}\frac{1}{z_1}\left(z_ka_{m_k} + \sum_{m\in \Delta\setminus S}m^kt^{w_m}z^ma_{m_k}\right)W^0_1(f).
\]

This computation together with lemma \ref{be} implies the desired result.
\hfill\rule{2.1mm}{2.1mm}\\
\begin{prop}
\label{bb}
Let $f$ be a bounded smooth function of terms in (\ref{bf}). Near $x\in X_t$ satisfying $S_x=S$, $f$ and its multi-derivatives with respect to $\{W_j,\bar{W}_j\}_{j=2}^l$ are bounded (uniformly with respect to $t$).
\end{prop}
{\bf Proof:} Since $\eta_m\leq 1$, by taking the constant $\kappa>0$ large, we can ensure that $a_m \geq C>0$. When $S_x=S$, lemma \ref{bd} gives us the filtrations $\{S_k\}_{k=1}^K$ and $\{M_k\}_{k=1}^K$. For $m\in \Delta\setminus S$ and $j\in I_m$, there exists $k$ such that $m\in M_k\setminus M_{k-1}$. Since $m^j\not=0$, we have $m_j\in S_k$. Lemma \ref{bd} then implies that $\eta_m \leq t_k^{b_k} \leq \eta_{m_j}^{b_k/c_k}$. Hence there exists a $C>0$ such that $a_{m_j} \leq Ca_m$. Consequently, the terms in (\ref{bf}) are bounded (uniformly with respect to $t$). By lemma \ref{bh}, we get the desired results.
\hfill\rule{2.1mm}{2.1mm}\\

\se{Construction and estimates of the approximate metrics}
Choose the Fubini-Study metric $\omega_{FS}$ on $\mathbb{CP}^{|\Delta|-1}$, and define

\[
\hat{\omega}_t = e_t^*\omega_{FS}.
\]

Then the \k potential of $\hat{\omega}_t$ is the logarithm of the volume form

\[
\hat{V}_t = \sum_{m\in \Delta} \Omega_{t,m}\otimes\bar{\Omega}_{t,m}.
\]

Let

\[
h_t = \kappa^2\sum_{S\in Z_{\rm top}}A_S^2,\ A_S = a_S^{-1},\ a_S = \prod_{m\in S} a_m,
\]
then

\[
\partial A_S = - A_S \partial \log a_S,\ \ \partial \bar{\partial} A_S = - A_S \partial \bar{\partial} \log a_S + A_S \partial \log a_S\bar{\partial} \log a_S,
\]
\[
\frac{i}{2\pi}\partial \bar{\partial} \log h_t = -\frac{i}{2\pi}\sum_{S\in Z_{\rm top}} \lambda_S \partial \bar{\partial} \log a_S^2 + \omega' = \sum_{S\in Z_{\rm top}} \lambda_S \omega_S + \left(\sum_{m\in \Delta} \frac{2\lambda_m}{a_m}\right)\hat{\omega}_t + \omega', 
\]

where

\[
\omega_S = \frac{i}{\pi}\sum_{m\in S}\frac{\partial a_m \bar{\partial} a_m}{a_m^2},\ \ \lambda_S = A_S^2\left(\sum_{S\in Z_{\rm top}} A_S^2\right)^{-1},\ \ \lambda_m = \sum_{S\ni m} \lambda_S,
\]
\[
\omega' = \frac{i}{2\pi}\sum_{S\in Z_{\rm top}} \lambda_S \partial \log a_S^2 \bar{\partial} \log a_S^2 - \frac{i}{2\pi}\left(\sum_{S\in Z_{\rm top}} \lambda_S \partial \log a_S^2\right)\left(\sum_{S\in Z_{\rm top}} \lambda_S \bar{\partial} \log a_S^2\right) 
\]
\[
= \frac{i}{4\pi}\sum_{S,S'\in Z_{\rm top}} \lambda_S\lambda_{S'} \partial (\log a_S^2 - \log a_{S'}^2) \bar{\partial} (\log a_S^2 - \log a_{S'}^2).
\]

Let $V_t = h_t\hat{V}_t$, then we have the \k form of the approximate metric

\[
\omega_t = \frac{i}{2\pi}\partial\bar{\partial}\log V_t = \hat{\omega}_t + \frac{i}{2\pi}\partial\bar{\partial}\log h_t
\]
\[
= \left(1 + \sum_{m\in \Delta} \frac{2\lambda_m}{a_m}\right)\hat{\omega}_t + \sum_{S\in Z_{\rm top}} \lambda_S \omega_S + \omega'.
\]

\begin{lm}
\label{cg}
Near $x\in X_t$ satisfying $S_x=S$, $\omega_t$ is quasi-isometric to
\[
\omega^\circ_S = \frac{i}{\pi}\sum_{i=2}^l\frac{1}{a_{m_i}^2}\frac{dz_id\bar{z}_i}{|z_i|^2}.
\]
\end{lm}
{\bf Proof:}
\[
\partial a_{m_i} = \frac{dz_i}{z_i} + \sum_{k=1}^l \left(\|s_{m_k}\|_t^2 + \sum_{m\in \Delta\setminus S} m^k\|s_m\|_t^2\right)\frac{dz_k}{z_k}
\]

Since $\|s_m\|_t^2 = O(t^a)$ for $m\in \Delta\setminus S$, we have

\[
\bigwedge_{i=1}^l\partial a_{m_i} = \left(1+ \sum_{k=1}^l \|s_{m_k}\|_t^2 + O(t^a)\right)\bigwedge_{i=1}^l\frac{dz_i}{z_i}.
\]

Consequently, $\displaystyle \frac{i}{\pi}\sum_{i=1}^l\frac{\partial a_{m_i} \bar{\partial}a_{m_i}}{a_{m_i}^2}$ is quasi-isometric to $\displaystyle \frac{i}{\pi}\sum_{i=1}^l\frac{1}{a_{m_i}^2}\frac{dz_id\bar{z}_i}{|z_i|^2}$. It is straightforward to verify the following:\\

(1) $\displaystyle \frac{i}{\pi}\sum_{i=1}^l\frac{\partial a_{m_i} \bar{\partial}a_{m_i}}{a_{m_i}^2}$ is quasi-isometric to $\displaystyle \sum_{S\in Z_{\rm top}} \lambda_S \omega_S + \omega'$,\\

(2) $\displaystyle \frac{i}{\pi}\sum_{i=1}^l\frac{1}{a_{m_i}^2}\frac{dz_id\bar{z}_i}{|z_i|^2}$ dominates $\hat{\omega}_t$,\\

(3) $\displaystyle \frac{i}{\pi}\sum_{i=1}^l\frac{1}{a_{m_i}^2}\frac{dz_id\bar{z}_i}{|z_i|^2}$ is quasi-isometric to $\displaystyle \frac{i}{\pi}\sum_{i=2}^l\frac{1}{a_{m_i}^2}\frac{dz_id\bar{z}_i}{|z_i|^2}$.\\

Combining these, we get the desired conclusion.
\hfill\rule{2.1mm}{2.1mm}\\
\begin{prop}
\label{ck}
Near $x\in X_t$ satisfying $S_x=S$, $\omega_t$ is $C^{\infty}$-quasi-isometric to $\omega^\circ_S$, and the basis $\{W_j,\bar{W}_j\}_{j=2}^l$ is proper with respect to both metrics.
\end{prop}
{\bf Proof:} It is straightforward to check that the coefficients of $[W_j,W_k]$, $[W_j,\bar{W}_k]$, $[\bar{W}_j,\bar{W}_k]$ with respect to the basis $\{W_j,\bar{W}_j\}_{j=2}^l$ are all bounded smooth functions of terms in (\ref{bf}). For example,

\[
[W_j,W_k] = W_j(\log a_{m_k})W_k - W_k(\log a_{m_j})W_j.
\]
\[
W_j(\log a_{m_k}) = W^0_j(\log a_{m_k}) -\frac{a_{m_j}}{a_{m_1}z_1}\left(z_j + \sum_{m\in \Delta\setminus S}m^jt^{w_m}z^m\right)W^0_1(\log a_{m_k}).
\]

By (\ref{ba}), $W_j(\log a_{m_k})$ is clearly a bounded smooth function of terms in (\ref{bf}). By proposition \ref{bb} and proposition \ref{ga}, we have that $\{W_j,\bar{W}_j\}_{j=2}^l$ is proper with respect to $\omega^\circ_S$.\\

It is straightforward to check that the coefficients of $\omega_t$ with respect to the basis $\{W_j,\bar{W}_j\}_{j=2}^l$ are all bounded smooth functions of terms in (\ref{bf}). By proposition \ref{bb}, lemma \ref{cg} and proposition \ref{ga}, we have that $\{W_j,\bar{W}_j\}_{j=2}^l$ is also proper with respect to $\omega_t$. Consequently, $\omega_t$ is $C^{\infty}$-quasi-isometric to $\omega^\circ_S$.
\hfill\rule{2.1mm}{2.1mm}\\
\begin{prop}
\label{cc}
The curvature of $g_t$ and its derivatives are all uniformly bounded with respect to $t$.
\end{prop}
{\bf Proof:} This is a direct consequence of propositions \ref{ck} and \ref{ga}.
\hfill\rule{2.1mm}{2.1mm}\\
\begin{lm}
\label{cf}
Near $x\in X_t$ satisfying $S_x=S$, 
\[
\hat{V}_t = (1+|z|^2 + O(t^a))\Omega_{t,m_0}\otimes\bar{\Omega}_{t,m_0} = (1 + O(t^a))\frac{1+|z|^2}{|z_1|^2} \prod_{i=2}^l \frac{dz_id\bar{z}_i}{|z_i|^2}.
\]
\end{lm}
{\bf Proof:} Notice that $t^{w_m}z^m = O(t^a)$ for $m\in \Delta\setminus S$. The first equality is obvious from the formula

\[
\hat{V}_t = (1+|z|^2 + \sum_{m\in \Delta\setminus S} |t^{w_m}z^m|^2)\Omega_{t,m_0}\otimes\bar{\Omega}_{t,m_0}.
\]

Since $\displaystyle s_t/s_{m_0} = 1+ z_1 + \cdots + z_l + \sum_{m\in \Delta\setminus S} t^{w_m}z^m$, $\displaystyle d(s_t/s_{m_0}) = \sum_{i=1}^l (z_i + O(t^a))\frac{dz_i}{z_i}$. The second equality is a direct consequence of

\[
\Omega_{t,m_0} = (\Omega \otimes (d(s_t/s_{m_0}))^{-1})|_{X_t} = (1 + O(t^a))\frac{1}{z_1} \prod_{i=2}^l \frac{dz_i}{z_i}.
\]
\hfill\rule{2.1mm}{2.1mm}\\

Assume $\displaystyle e^{-\phi_t} = \frac{\omega_t^{l-1}}{V_t}$, we have\\
\begin{prop}
\label{cb}
$|\phi_t|$ is bounded independent of $t$.
\end{prop}
{\bf Proof:} According to lemmas \ref{cf} and \ref{cg} and the definition of $h_t$, we have

\[
e^{-\phi_t} = \frac{\omega_t^{l-1}}{V_t} = \frac{\omega_t^{l-1}}{h_t\hat{V}_t} \sim 1.
\]

Therefore $|\phi_t|$ is bounded independent of $t$. 
\hfill\rule{2.1mm}{2.1mm}\\
\begin{prop}
\label{cd}
For any $k$, $\|\phi_t\|_{C^k,g_t}$ is uniformly bounded with respect to $t$.\\
\end{prop}
{\bf Proof:} With proposition \ref{cb}, it is straightforward to check that $\phi_t$ is a bounded smooth function of terms in (\ref{bf}). Applying proposition \ref{bb} and proposition \ref{ck}, we get the desired result.
\hfill\rule{2.1mm}{2.1mm}\\

Let $\psi_t: \tilde{X}_0\setminus{\rm Sing}(\tilde{X}_0) \rightarrow X_t$ be the Hamiltonian-gradient flow (lifted to $\tilde{X}_0\setminus{\rm Sing}(\tilde{X}_0)$ from $X_0\setminus{\rm Sing}(X_0)$) with respect to the family $X_t \subset \mathbb{CP}^{|\Delta|-1}$ under the Fubini-Study metric. Then $\psi_0 = e_0$. As discussed in \cite{sl1} and the references therein, the Hamiltonian-gradient flow $\psi_t$ for a general family $\{X_t,g_t\}$ with total space $({\cal X},g)$ satisfying $g|_{X_t} = g_t$ is defined along radial direction $I_{\theta_0} = \{t=re^{i\theta_0}:0\leq r<1\}$, and is completely determined by the restriction $\{X_t,g_t\}_{t\in I_{\theta_0}}$. The flow $\psi_t$ for a similar situation was discussed in the remark at the end of \cite{ke1}. The fact that we will need of the flow is that locally near any compact set $F\in \tilde{X}_0\setminus{\rm Sing}(\tilde{X}_0)$, the flow $\{\psi_t\}_{t\in I_{\theta_0}}$ is a bounded smooth family of smooth diffeomorphisms, which is obvious because the Hamiltonian-gradient vector field is smooth in such region.\\

Let $h_0|_{\tilde{X}_{0,S}} = h_{0,S}|_{\tilde{X}_{0,S}} = \kappa^2A_S^2$ for each $S\in Z_{\rm top}$, and $\tilde{\omega}_0 = \hat{\omega}_0 + \frac{i}{2\pi}\partial\bar{\partial}\log h_0$, where $\hat{\omega}_0 = e_0^*\omega_{FS}$. Then we have\\
\begin{prop}
\label{dd}
The approximate metric $g_t$ on $X_t$ will converge to the complete metric $\tilde{g}_0$ on $\tilde{X}_0\setminus{\rm Sing}(\tilde{X}_0)$ in the sense of Cheeger-Gromov: for any compact subset $F \subset \tilde{X}_0\setminus{\rm Sing}(\tilde{X}_0)$, $\psi_t^*g_t$ converge to $\tilde{g}_0$ on $F$ in $C^k$-topology on the space of Riemannian metrics as $t$ goes to $0$ for any $k$.
\end{prop}
{\bf Proof:} Without loss of generality, we may assume that $F$ is a compact subset in $\tilde{X}_{0,S}$ for fixed $S\in Z_{\rm top}$. Since $e_0(F)$ is a compact subset in $T_S$, there exists a constant $C>0$ and a choice of small neighborhood $U_S \subset \mathbb{CP}^{|\Delta|-1}$ of $e_0(F)$ such that $1 \geq|z_1|\geq\cdots \geq |z_l| \geq C$ in $U_S$. Hence, for $x\in X_t\cap U_S$, we have $S_x = S$. And it is straightforward to see that $\{W_j,\bar{W}_j\}_{j=2}^l$ is equivalent to $\{\frac{\partial}{\partial z_i},\frac{\partial}{\partial \bar{z}_i}\}_{i=2}^l$ in $X_t\cap U_S$. By propositions \ref{gc} and \ref{gb}, we have that $g_t$ is $C^\infty$-quasi-isometric to $\hat{g}_t$ in $X_t\cap U_S$, uniform with respect to $t$. It is also straightforward to see that $\{\psi_t^*\frac{\partial}{\partial z_i},\psi_t^*\frac{\partial}{\partial \bar{z}_i}\}_{i=2}^l$ are equivalent to each other for different $t$. Hence $\psi_t^*\hat{g}_t$ are $C^\infty$-quasi-isometric to each other for different $t$. Consequently, $\psi_t^*g_t$ are $C^\infty$-quasi-isometric to each other uniformly for different $t$. It is straightforward to see that $\displaystyle\lim_{t\rightarrow 0}\psi_t^*g_t = \tilde{g}_0$ on $F$ as tensors under $C^0$-topology. By proposition \ref{gd}, we have $\psi_t^*g_t$ converge to $\tilde{g}_0$ on $F$ in $C^k$-topology on the space of Riemannian metrics as $t$ goes to $0$ for any $k$.
\hfill\rule{2.1mm}{2.1mm}\\

\se{Construction and degeneration of the \ke metrics}
In \cite{Tian1}, using the Monge-Amp\`{e}re estimate of Aubin and Yau, Tian essentially proved the following theorem when $X_t$ are compact. The case when $X_t$ are complete can be proved in exactly the same way if one uses the Cheng-Yau's Monge-Amp\`{e}re estimate for complete manifold instead (\cite{CY1}).\\
\begin{theorem}
\label{db}
Assume that $\phi_t$, the curvature of $g_t$ and their multi-derivatives are all bounded uniformly independent of $t$, the approximate metrics $g_t$ on $X_t$ converge to the complete metric $\tilde{g}_0$ on $\tilde{X}_0\setminus{\rm Sing}(\tilde{X}_0)$ in the sense of Cheeger-Gromov, then the complete \ke metric $g^{\rm KE}_t$ on $X_t$ will converge to the complete \ke metric $g_{E,0}$ on $\tilde{X}_0\setminus{\rm Sing}(\tilde{X}_0)$ in the sense of Cheeger-Gromov.
\hfill\rule{2.1mm}{2.1mm}\\
\end{theorem}
\begin{theorem}
\label{dc}
The complete Cheng-Yau \ke metrics $g^{\rm KE}_t$ on $X_t$ satisfying ${\rm Ric}(g^{\rm KE}_t) = - g^{\rm KE}_t$ converge to the complete Cheng-Yau \ke metric $g_{E,0}$ on $\tilde{X}_0\setminus{\rm Sing}(\tilde{X}_0)$ in the sense of Cheeger-Gromov.
\end{theorem}
{\bf Proof:} This theorem is a direct corollary of theorem \ref{db} and propositions \ref{cc}, \ref{cd}, \ref{dd}. 
\hfill\rule{2.1mm}{2.1mm}\\

It is easy to see that our construction actually implies the following asymptotic description of the family of \ke metrics.\\
\begin{theorem}
\label{da}
\ke metric $g^{\rm KE}_t$ on $X_t$ is uniformly $C^\infty$-quasi-isometric to the explicit approximate metric $g_t$. More precisely, there exist constants $C_1,C_2>0$ independent of $t$ such that $C_1 g_t \leq g^{\rm KE}_t\leq C_2 g_t$ and $g^{\rm KE}_t$ is uniformly $C^\infty$-bounded with respect to $g_t$.
\end{theorem}
{\bf Proof:} 
The uniform $C^0$-estimate of the complex Monge-Amp\`{e}re equations implies that $C_1 \omega_t^{l-1} \leq (\omega^{\rm KE}_t)^{l-1}\leq C_2 \omega_t^{l-1}$ for some $C_1,C_2>0$. The uniform $C^2$-estimate of the complex Monge-Amp\`{e}re equations implies that ${\rm Tr}_{g_t} g^{\rm KE}_t$ is uniformly bounded from above. Consequently, $g^{\rm KE}_t$ is uniformly quasi-isometric to $g_t$. The uniform complex Monge-Amp\`{e}re estimates based on the uniform estimates in propositions \ref{cb}, \ref{cc}, \ref{cd} imply that $g^{\rm KE}_t$ is uniformly $C^\infty$-bounded with respect to $g_t$.
\hfill\rule{2.1mm}{2.1mm}\\

\se{Toroidal degeneration}
Although the degenerations discussed in this paper are examples of toroidal degenerations, due to the special structure of our situation, we are able to construct the degeneration family of \ke metrics directly without incurring the general theory of toroidal degeneration. Therefore, the discussion in this section is not absolutely necessary for other sections of this paper, but will serve to reveal the underlying algebraic degeneration structure.\\

As one may observe, the degenerations discussed in our paper are somewhat more general than usual algebraic degenerations. In algebraic geometry, degenerations are usually over a disc $D = \{|t|< 1\} \subset \mathbb{C}$ with the central fibre $X_0$ over $t=0$ being the singular fibre. In our situation, since $w$ is not necessarily rational, there are multi-fibres $X_t$ over each $t\not=0$, while there is only one singular fibre $X_0$ over $t=0$. We will call such degeneration $\mathbb{R}$-degeneration. To avoid multi-fibres, we may consider the restriction of the degeneration to a ray $I_{\theta_0}=\{t=re^{i\theta_0}:0\leq r <1\}$ for a fixed angle $\theta_0$ or a fan-like region $D_{\theta_0,\epsilon} = \{t=re^{i\theta}:0\leq r <1,\theta_0-\epsilon \leq \theta \leq \theta_0+\epsilon\}$.\\ 

Recall that $M$ is an integral lattice, and $N = M^\vee$. Consider $w = \{w_m\}_{m\in \Sigma''(1)}$, where $\Sigma''(1)$ is a finite subset of primitive elements in $M$ whose real convex span is $M_{\mathbb{R}}$, and $w$ is convex at the origin (in another word, $w$ can be adjusted by a linear function on $M$ such that $w_m>0$ for $m\in \Sigma''(1)$). Then there exists the maximal piecewise linear function $p$ on $M$ satisfying $p(m)\leq w_m$ for $m\in \Sigma''(1)$. Let $\Sigma'(1)$ be the subset of $\Sigma''(1)$ such that $p(m)= w_m$, and $\Sigma(1)$ be the subset of $\Sigma'(1)$ such that $p$ is not linear near $m$ for $m\in \Sigma'(1)$. $p$ determines a complete fan $\Sigma$ on $M$ with $\Sigma(1)$ identified with the set of 1-dimensional cones in $\Sigma$. We will assume $w$ to be generic, then $\Sigma$ will be a simplicial fan. $w$ is called convex (resp. strictly convex) if $\Sigma''(1)=\Sigma'(1)$ (resp. $\Sigma''(1)=\Sigma(1)$). In particular, $w' = w|_{\Sigma'(1)}$ is convex and $w'' = w|_{\Sigma(1)}$ is strictly convex.\\ 

For $t\not=0$, $i_t = \{s_m(z)\}_{m\in \Sigma''(1)}: N_{\mathbb{C}^*} \rightarrow \mathbb{C}^{|\Sigma''(1)|}$ with $s_m(z) = t^{w_m}z^m$ defines a family of toric embedding. Let $X_t=i_t(N_{\mathbb{C}^*})$ and $\displaystyle X_0 = \bigcup_{\sigma\in \Sigma} X_{0,\sigma}\subset \mathbb{C}^{|\Sigma'(1)|}$, where $X_{0,\sigma} = i_{0,\sigma}(N_{\mathbb{C}^*})$ with $i_{0,\sigma} = \{z^m\}_{m\in \Sigma'(1)\cap\sigma}: N_{\mathbb{C}^*} \rightarrow \mathbb{C}^{|\Sigma'(1)\cap\sigma|} \subset \mathbb{C}^{|\Sigma''(1)|}$. $X_0$ is the set-theoretical limit of $X_t$ in $\mathbb{C}^{|\Sigma''(1)|}$ as $t\rightarrow 0$. Such degeneration $\{X_t\}$ is called toric degeneration (determined by $w$). Two toric degenerations with the same $X_0$ are said to be equivalent to each other. Among all the toric degenerations with the same $w'' = w|_{\Sigma(1)}$, there is a minimal equivalent class determined by toric degeneration satisfying $w'' = w$, there is also a maximal equivalent class determined by toric degenerations satisfying ${\rm Span}_{\mathbb{Z}}(\Sigma'(1)\cap \sigma) = {\rm Span}_{\mathbb{Z}}(\sigma)$ for all $\sigma \in \Sigma$. General toric degenerations degenerate $X_t$ to different components of $X_0$ with possibly different multiplicities, while the maximal toric degeneration is simple, namely, it degenerate $X_t$ to $X_0$ with multiplicity 1.\\
 
A degeneration $\{X_t\}$ is called toroidal if it is locally a toric degeneration times a smooth manifold.\\
\begin{lm}
$\{X_t\}$ and $\{Y_t\}$ defined in section 3 are toroidal degenerations at $t=0$.
\end{lm}
{\bf Proof:} Since the hypersurface $H$ intersects each smooth strata of $Y_0$ transversely, we only need to verify the lemma for $\{Y_t\}$.\\

For $S \in Z$, $w$ defined in the beginning of section 3 can be adjusted by affine function on $M$ so that $w_m=0$ for $m\in S$ and $w_m >0$ for $m\in \Delta\setminus S$. Let $F$ be a compact subset in $T_S$. Then for $x\in Y_t$ near $F$, we have $S\subset S_x = \{m_0,m_1,\cdots,m_l\}$. Assume $S = \{m_0,m_1,\cdots,m_{l'}\}$ for some $l'\leq l$. Then locally near $F$, we have the map $i_t = (i'_t,i''_t): Y_t \rightarrow \mathbb{C}^{|\Delta|-1}$ with $i'_t = (z_1,\cdots,z_{l'})$ and $i''_t = \{t^{w_m}z^m\}_{m\in \Delta\setminus S}$. Since $F\subset T_S$ is compact, we have $C_1\geq |z_i| \geq C_2>0$ for $1\leq i \leq l'$. Consequently, $i''_t$ is naturally a toric degeneration and $\{Y_t\}$ near $F$ is a toric degeneration family times $i_{0,S}^{-1}(F)$, where $i_{0,S}: (N_S)_{\mathbb{C}^*} \rightarrow T_S$ is a finite cover.
\hfill\rule{2.1mm}{2.1mm}\\

It is not hard to check that the canonical limit $\tilde{X}_0$ (resp. $Y_Z$) defined in section 3, when restricted to each local toric degeneration model, is exactly the simple maximal toric degeneration limit of $\{X_t\}$ (resp. $\{Y_t\}$).\\

\se{The minimal Lagrangian vanishing torus}
In \cite{hl}, we constructed minimal Lagrangian vanishing torus for toroidal degeneration family of \ke manifolds discussed in \cite{ke2} near maximal degeneracy points. In this section, we will apply results in \cite{hl} to construct minimal Lagrangian vanishing torus in $(X_t,g^{\rm KE}_t)$ near maximal degeneracy points in $X_0$. The maximal degeneracy points in $X_0$ are $0$-dimensional stratas in $X_0$. They are the intersections of $1$-dimensional stratas in $Y_0$ and $H$. $1$-dimensional stratas in $Y_0$ correspond to $1$-simplices in $Z$. Let $S_1 = \{m_0,m_1\} \in Z$ be a $1$-simplex. Let $\hat{M} = M/\mathbb{Z}(m_1-m_0)$, with the projection $\pi: M \rightarrow \hat{M}$. We will use $\hat{\Sigma}(1) \subset \hat{M}$ to denote the image of $\Delta_{S_1} = \displaystyle \bigcup_{S_1\subset S \in Z} (S\setminus S_1) \subset M$ into $\hat{M}$. Let $\hat{Z}$ denote the set of simplices $\pi(S\setminus S_1)$ in $\hat{M}$ for $S\in Z$ satisfying $S_1 \subset S$. Adjust $\{w_m\}_{m\in \Delta}$ by affine function, we may assume $w_{m_0} = w_{m_1}=0$ and $w_m\geq \tilde{w}$ for $m\in \Delta\setminus S_1$ and a constant $\tilde{w}>0$.\\

For $m\in \Delta_{S_1}$, define $\hat{s}_{\hat{m}} = (s_m/s_{m_0})|_{\hat{T}}$, where $\hat{m} = \pi(m)$, $\hat{T} = \{x\in (\mathbb{C}^*)^l: s_{m_1}(x)= -s_{m_0}(x)\}$. $\{t^{w_m}\hat{s}_{\hat{m}}\}_{\hat{m}\in \hat{\Sigma}(1)}: \hat{T} \rightarrow \mathbb{C}^{|\hat{\Sigma}(1)|}$ defines a toric degeneration. Let $\hat{N} = \hat{M}^{\vee} = (m_1-m_0)^\perp$. There is the natural moment map $F: \hat{T} \cong \hat{N}_{\mathbb{C}^*}\rightarrow \hat{N}_{\mathbb{R}}$. Recall from \cite{hl}

\[
\Delta_\tau = \{x \in \hat{N}_{\mathbb{R}}|\langle \hat{m},x\rangle + \tau w_m \geq 0,\ {\rm for}\ \hat{m}\in \hat{\Sigma}(1)\},\ \ {\rm where}\ \tau = -\log |t|.
\]

On $F^{-1}(\Delta_{\tau}) \subset \hat{T}$, define a family of toric \k metrics $g^{\rm tor}_t$ with \k potential $\rho_\tau = \log \hat{h}_t$, where
\[
\hat{h}_t = \sum_{S\in \hat{Z}_{\rm top}}\hat{A}_S^2,\ \hat{A}_S = \hat{a}_S^{-1},\ \hat{a}_S = \prod_{\hat{m}\in S} \hat{a}_{\hat{m}},\ \hat{a}_{\hat{m}} = -\log |t^{w_m}\hat{s}_{\hat{m}}|^2.
\]

Let $O$ be a maximal degeneracy point in $X_0$. Then locally near $O$, $X_t$ can be identified with $F^{-1}(\Delta_{\tau})$. More precisely, we will use the coordinate $\{z_k = s_{m_k}/s_{m_0}\}_{k=2}^l$ with respect to a fixed $S_* = \{m_0,\cdots,m_l\}(\cong \{m_0,m_1\}\cup\hat{S}_*)\in Z_{\rm top}$ as toric coordinate for the identification, where $\hat{S}_*\in \hat{Z}_{\rm top}$.\\
\begin{lm}
\label{fd}
There exists constant $\mu>0$ so that $g^{\rm tor}_t$ is $C^\infty$-quasi-isometric to $g_t$ on $F^{-1}(\Delta_{\tau-\mu})$ uniformly for $\tau$ large.
\end{lm}
{\bf Proof:} For $x\in X_t \subset N_{\mathbb{C}^*}$, Let $S_x\in Z_{\rm top}$ denote the simplex chosen in lemma \ref{bd}. When $x$ is near $O$, $S_x$ contains $S_1 = \{m_0,m_1\}$. Assume that $S_x = \{m_0,\cdots,m_l\} \cong \{m_0,m_1\}\cup\hat{S}_x$, where $\hat{S}_x\in \hat{Z}_{\rm top}$. Choose $\{z_k=s_{m_k}/s_{m_0}\}_{k=1}^l$ and $z = \{z_k\}_{k=2}^l$ to be coordinates of $N_{\mathbb{C}^*}$ and $X_t$ near $O$. Choose $\hat{z} = \{\hat{z}_k = \hat{s}_{\hat{m}_k}\}_{k=2}^l$ to be the coordinate for $F^{-1}(\Delta_{\tau}) \subset \hat{T}$. Recall that $X_t$ near $O$ is identified with $F^{-1}(\Delta_{\tau})$ using the fixed $S_*$. $z$ and $\hat{z}$ are identified under such identification if $S_x=S_*$. In general, $S_x$ is not necessarily equal to $S_*$ and the coordinate transformation is generally, $\hat{z}_j = z_j (-z_1)^{d_j}$ for $2\leq j \leq l$. $z_j\frac{\partial}{\partial z_j} = \hat{z}_j\frac{\partial}{\partial\hat{z}_j} + d_k\frac{z_j}{z_1}\frac{\partial z_1}{\partial z_j}\hat{z}_k\frac{\partial}{\partial\hat{z}_k}$.\\

By proposition \ref{ck}, $\{W_j,\bar{W}_j\}_{j=2}^l$ is a proper basis for $g_t$. With similar (and slightly simpler) arguments as in the proof of proposition \ref{ck}, one can see that $\{\hat{W}_j,\bar{\hat{W}}_j\}_{j=2}^l$ is a proper basis for $g^{\rm tor}_t$, where $\hat{W}_j = \hat{a}_{\hat{m}_j} \hat{z}_j\frac{\partial}{\partial\hat{z}_j}$. We have

\[
W_j = b_j^k \hat{W}_k = \frac{a_{m_j}}{\hat{a}_{\hat{m}_j}}\hat{W}_j + d_k\frac{a_{m_j}z_j}{\hat{a}_{\hat{m}_k}z_1}\frac{\partial z_1}{\partial z_j}\hat{W}_k.
\]
\[
\det(b_j^k) = \frac{a_{\hat{S}_x}}{\hat{a}_{\hat{S}_x}}\left(1 + \sum_{j=2}^ld_j\frac{z_j}{z_1}\frac{\partial z_1}{\partial z_j}\right),\ \ {\rm where}\ a_{\hat{S}_x} = \prod_{m\in S_x\setminus S_1} a_m.
\]

$\hat{z}\in F^{-1}(\Delta_{\tau-\mu})$ implies that $\hat{a}_m \geq \tilde{w}\mu$ and $|t^{w_m}\hat{z}^{\hat{m}}|\leq e^{-\tilde{w}\mu}$ for $m \in \Delta_{S_1}$. Since $\log |z_1|^2$ is bounded by lemma \ref{bg}, we also have $|t^{w_m}z^m|= O(e^{-\tilde{w}\mu})$ for $m \in \Delta\setminus S_1$. Notice that
\[
a_{m_0} = \kappa + \log\left(1 + |z_1|^2 + \sum_{m\in \Delta\setminus S_1}|t|^{2w_m}|z^m|^2\right)
\]
is also bounded. It is easy to verify that

\[
\frac{a_{m_j}}{\hat{a}_{m_j}} = 1 + \frac{a_{m_0} -d_j\log |z_1|^2}{\hat{a}_{m_j}}.
\]

When $\mu$ is large enough, we have $a_{m_j}/\hat{a}_{m_j} = 1 + O(1/\mu)$ for $2\leq j \leq l$. Consequently, $a_{\hat{S}_x}/\hat{a}_{\hat{S}_x} = 1 + O(1/\mu)$. On the other hand,
\[
z_j\frac{\partial z_1}{\partial z_j} = t^{w_{m_j}}z_j + \sum_{m\in \Delta\setminus S_x}m^jt^{w_m}z^m = O(e^{-\tilde{w}\mu}).
\]

Consequently, when $\mu$ is large, we have $\det(b_j^k) = 1 + O(1/\mu,e^{-\tilde{w}\mu})\geq C>0$.\\

It is straightforward to verify that $b_j^k$ are bounded smooth functions of terms in (\ref{bf}) using above computations and proposition \ref{aa}. Then proposition \ref{bb} implies that $b_j^k$ and their multi-derivatives with respect to $\{W_j,\bar{W}_j\}_{j=2}^l$ are bounded. Consequently, $\{W_j,\bar{W}_j\}_{j=2}^l$ is equivalent to $\{\hat{W}_j,\bar{\hat{W}}_j\}_{j=2}^l$. By propositions \ref{gc} and \ref{gb}, we get the desired conclusion.
\hfill\rule{2.1mm}{2.1mm}\\

Let $\{g^{\rm tor}_t\}$ be a smooth family of toric \k metrics and $\{g^{\rm KE}_t\}$ be a smooth family of \ke metrics defined on the family of spaces $\{F^{-1}(\Delta_{\tau})\}$. Theorem 4.1 in \cite{hl} can be rephased as the following:\\
\begin{theorem}
\label{fe}
For certain fixed $\mu$, if $g^{\rm tor}_t$ and $g^{\rm KE}_t$ defined on $F^{-1}(\Delta_{\tau})$ are $C^\infty$-quasi-isometric on $F^{-1}(\Delta_{\tau-\mu})$ (uniform with respect to $t$), and the \k potential $\rho_\tau$ of the toric \k metric $g^{\rm tor}_t$ (viewed as a function on $\Delta_\tau$) satisfies $\rho_\tau(x) - \rho_1(x/\tau) = C(\tau)$ for $x\in \Delta_\tau$ and $\displaystyle \lim_{c\rightarrow 1} \rho_1|_{\partial\Delta_c} = +\infty$, then there exists a smooth family of minimal Lagrangian torus $L_t \subset (X_t,\omega^{\rm KE}_t)$.
\hfill\rule{2.1mm}{2.1mm}
\end{theorem} 
\begin{theorem}
\label{fa}
Let $O$ be a maximal degeneracy point in $X_0$. Then there exists a smooth family of minimal Lagrangian torus $L_t\subset (X_t,\omega^{\rm KE}_t)$ for $t$ small that approaches $O$ when $t$ approaches $0$.
\end{theorem}
{\bf Proof:} We will take $g^{\rm tor}_t$ as the toric metric defined earlier on $F^{-1}(\Delta_{\tau})$. As discussed earlier in this section, $X_t$ near $O$ can be identified with $F^{-1}(\Delta_{\tau})$. We will take $g^{\rm KE}_t$ to be the restriction to $F^{-1}(\Delta_{\tau})$ of the \ke metric on $X_t$. Lemma \ref{fd} implies that $g^{\rm tor}_t$ is $C^\infty$-quasi-isometric to $g_t$ on $F^{-1}(\Delta_{\tau-\mu})$. Theorem \ref{da} implies that $g_t$ is $C^\infty$-quasi-isometric to $g^{\rm KE}_t$. Consequently, $g^{\rm tor}_t$ is $C^\infty$-quasi-isometric to $g^{\rm KE}_t$ by proposition \ref{ge}. From the definition of $\rho_\tau$, it is easy to observe that $\rho_\tau(x) = \rho_1(x/\tau) - 2(l+1)\log \tau$ and $\displaystyle \lim_{c\rightarrow 1} \rho_1|_{\partial\Delta_c} = +\infty$, apply theorem \ref{fe}, we get the desired result.
\hfill\rule{2.1mm}{2.1mm}\\

\ifx\undefined\bysame
\newcommand{\bysame}{\leavevmode\hbox to3em{\hrulefill}\,}
\fi

\end{document}